\documentclass[12pt,reqno]{article}
\usepackage{amsmath, amsthm, amscd, amsfonts, amssymb, graphicx, color}
\usepackage[bookmarksnumbered, colorlinks, plainpages]{hyperref}
\usepackage{secdot}
\newtheorem{theorem}{Theorem}[section]
\newtheorem{lemma}[theorem]{Lemma}

\theoremstyle{definition}
\newtheorem{definition}[theorem]{Definition}

\theoremstyle{remark}

\numberwithin{equation}{section}

\newcommand{\be}{\begin{equation}}
\newcommand{\ee}{\end{equation}}
\newcommand{\bea}{\begin{eqnarray}}
\newcommand{\eea}{\end{eqnarray}}
\begin{document}

\begin{center}
{\Large{\textbf{Almost Kenmotsu Manifolds Admitting Certain Critical Metric}}}

\end{center}
\vspace{0.1 cm}

\begin{center}

Dibakar Dey \\

Department of Pure Mathematics,\\
University of Calcutta,
35 Ballygunge Circular Road,\\
Kolkata - 700019, West Bengal, India,\\
E-mail: deydibakar3@gmail.com \\
\end{center}

\vspace{0.3 cm}
\textbf{Abstract:} In the present paper, we introduce the notion of $\ast$-Miao-Tam critical equation on almost contact metric manifolds and studied on a class of almost Kenmotsu manifold. It is shown that if the metric of a $(2n + 1)$-dimensional $(k,\mu)'$-almost Kenmotsu manifold $(M,g)$ satisfies the $\ast$-Miao-Tam critical equation, then the manifold $(M,g)$ is $\ast$-Ricci flat and locally isometric to the Riemannian product of a $(n + 1)$-dimensional manifold of constant sectional curvature $-4$ and a flat $n$-dimensional manifold. Finally, an illustrative example is presented to support the main theorem. \\



\section{ \textbf{Introduction}}
\vspace{0.3 cm}

\par In differential geometry and mathematical physics, finding Riemannian metrics on a compact orientable manifold $M$ that provides constant scalar curvature is a very classical problem. The critical points of the total scalar curvature functional plays an important role in this problem. The critical points of the total scalar curvature functional $\mathcal{S} : \mathcal{M} \rightarrow \mathbb{R}$ defined by
\bea
\nonumber \mathcal{S}(g) = \int_{\mathcal{M}} r_g dv_g 
\eea
on a compact orientable manifold $(M,g)$ restricted to the set of all Riemannian metrics $\mathcal{M}$ of unit volume are Einstein (see \cite{besse}, where $r_g$ is the scalar curvature and $dv_g$ the elementary volume form of $g$.\\

\par Let $(M,g)$ be a compact orientable Riemannian manifold of unit volume and $g^\ast$ be any symmetric bilinear form on $M$. Then the linearization of the scalar curvature $\mathcal{L}_g(g^\ast)$ is given by
\bea
\nonumber \mathcal{L}_g(g^\ast) = - \Delta_g(tr_g g^\ast) + div(div(g^\ast)) - g(g^\ast,S_g),
\eea
where $\Delta_g$ is the Laplacian of $g$ and $S_g$ is it's Ricci tensor. The formal $L^2$-adjoint $\mathcal{L}_g^\ast$ of the linearized scalar curvature $\mathcal{L}_g$ is defined as
\bea
\nonumber \mathcal{L}_g^\ast (\lambda) = - (\Delta_g \lambda)g + \nabla_g^2 \lambda - \lambda S_g,
\eea
where $(\nabla_g^2 \lambda)(X,Y) = Hess_g \lambda(X,Y) = g(\nabla_X D\lambda,Y)$ is the Hessian of the smooth function $\lambda$ on $M$ and $D$ is the gradient operator.\\

\begin{definition}
On a Riemannian manifold $(M,g)$ of dimension $n > 2$ if there exist a non-zero smooth function $\lambda$ such that $\mathcal{L}_g^\ast (\lambda) = g$, that is,
\bea
- (\Delta_g \lambda)g + \nabla_g^2 \lambda - \lambda S_g = g, \label{1.1}
\eea
then $g$ is said to satisfy the Miao-Tam critical equation.
\end{definition}
In the above definition, if the potential function $\lambda$ is a non-zero constant, then the solution of the Miao-Tam critical equation are Einstein metrics. In \cite{mt}, Miao-Tam proved that any Riemannian metric $g$ satisfying (\ref{1.1}) have constant scalar curvature. Further, they proved that any connected, compact Einstein manifold with smooth boundary satisfying (\ref{1.1}) is isometric to a geodesic ball. In \cite{patra}, 
Patra and Ghosh proved that a complete $K$-contact metric satisfying the Miao-Tam critical equation is isometric to a unit sphere. In 2017, Wang \cite{ywang} studied this notion on three dimensional almost Kenmotsu manifold and proved that if the metric of a three dimensional $(k,\mu)'$-almost Kenmotsu manifold satisfies the Miao-Tam critical equation, then the manifold is locally isometric to $\mathbb{H}^3(-1)$.\\

In 1959, Tachibana\cite{st} introduced the notion of $\ast$-Ricci tensor on almost Hermitian manifolds. Later in \cite{ht}, Hamada defined the $\ast$-Ricci tensor of real hypersurfaces in non-flat complex space form by
  \bea
  S^\ast(X,Y) = g(Q^\ast X,Y) = \frac{1}{2}(trace\{\phi \circ R(X,\phi Y)\}) \label{1.2}
  \eea
 for any vector fields $X,\; Y$ on $M$, where $Q^\ast$ is the $(1,1)$ $\ast$-Ricci operator. The $\ast$-scalar curvature is denoted by $r^\ast$ and is defined by $r^\ast = trace(Q^\ast)$.
\begin{definition}
An almost contact metric manifold $(M,g)$ is called $\ast$-Ricci flat if the $\ast$-Ricci tensor $S^\ast$ vanishes identically.
\end{definition}

\par Over the last decade, geometers and mathematical physicists developed several notions related to the  $\ast$-Ricci tensor. In 2016, the notion of $\ast$-Ricci soliton (\cite{new5}) was introduced. Later in 2019, the notion of $\ast$-critical point equation \cite{dey} was introduced. In this paper, we introduce the notion of $\ast$-Miao-Tam critical equation on an almost contact metric manifold as follows:
\begin{definition}
The metric of an almost contact metric manifold $(M,g)$ of dimension $n \geq 3$ is said to satisfy the $\ast$-Miao-Tam critical equation if there is a non-zero smooth function $\lambda : M \rightarrow \mathbb{R}$ satisfying
\bea
- (\Delta_g \lambda)g + \nabla_g^2 \lambda - \lambda S^\ast = g, \label{1.3}
\eea
where $S^\ast$ is the $(0,2)$ $\ast$-Ricci tensor of $g$, provided $S^\ast$ is symmetric.
\end{definition}

Note that, the $\ast$-Ricci tensor is not symmetric in general. In a $(k,\mu)'$-almost Kenmotsu manifold, the $\ast$-Ricci tensor $S^\ast$ is symmetric(given later) and hence the above definition is well defined in $(k,\mu)'$-almost Kenmotsu manifolds. Further, the $\ast$-Ricci tensor is also symmetric in different kind of almost contact metric manifolds. This gives us a wide range to study this notion in contact geometry.\\

\par The paper is organized as follows: In section 2, we give some basic concept of $(k,\mu)'$-almost Kenmotsu manifolds. Section 3 deals with $(k,\mu)'$-almost Kenmotsu manifolds satisfying $\ast$-Miao-Tam critical equation. In section 4, an example is presented to verify the result.\\

\section{\textbf{$(k,\mu)'$-almost Kenmotsu manifolds }}

An odd dimensional differentiable manifold $M$ is said to have an almost contact structure, if it admits a $(1,1)$ tensor field $\phi$, a characteristic vector field $\xi$ and a 1-form $\eta$ satisfying (\cite{bl}, \cite{bll}),
\be
\phi^{2}=-I+\eta\otimes\xi,\;\; \eta(\xi)=1, \label{2.1}
\ee
where $I$ denote the identity endomorphism. Here also $\phi\xi=0$ and $\eta\circ\phi=0$; both can be derived from (\ref{2.1}) easily. \\
If a manifold $M$ with an almost contact  structure admits a Riemannian metric $g$ such that
\be
\nonumber g(\phi X,\phi Y)=g(X,Y)-\eta(X)\eta(Y),
\ee
for any vector fields $X$, $Y$ on $M$, then $M$ is said to be an almost contact metric manifold. The fundamental 2-form $\Phi$ on an almost contact metric manifold is defined by $\Phi(X,Y)=g(X,\phi Y)$ for any $X$, $Y$ on $M$. Almost contact metric manifold such that $\eta$ is closed and $d\Phi=2\eta\wedge\Phi$ are  called almost Kenmotsu manifolds (\cite{dp}, \cite{pv}).
 
  Let us denote the distribution orthogonal to $\xi$ by $\mathcal{D}$ and defined by $\mathcal{D}=Ker(\eta)=Im(\phi)$. In an almost Kenmotsu manifold, since $\eta$ is closed, $\mathcal{D}$ is an integrable distribution. \\
 Let $M$ be a $(2n+1)$-dimensional almost Kenmotsu manifold. We denote by $h=\frac{1}{2}\pounds_{\xi}\phi$ and $l=R(\cdot, \xi)\xi$ on $M$. The tensor fields $l$ and $h$ are symmetric operators and satisfy the following relations \cite{pv}:
 \bea
h\xi=0,\;l\xi=0,\;tr(h)=0,\;tr(h\phi)=0,\;h\phi+\phi h=0, \label{2.2}
 \eea
\bea
 \nabla_{X}\xi=X - \eta(X)\xi - \phi hX(\Rightarrow \nabla_{\xi}\xi=0), \label{2.3}
\eea
 \be
  \phi l \phi-l = 2(h^{2} - \phi^{2}),\label{2.4}
 \ee
 \be
   R(X,Y)\xi = \eta(X)(Y - \phi hY) - \eta(Y)(X - \phi hX)+(\nabla_{Y}\phi h)X - (\nabla_{X}\phi h)Y, \label{2.5}
 \ee
 for any vector fields $X,Y$. The $(1,1)$-type symmetric tensor field $h'=h\circ\phi$ is anti-commuting with $\phi$ and $h'\xi=0$. Also it is clear that (\cite{dp}, \cite{waa})
 \bea
  h=0\Leftrightarrow h'=0,\;\;h'^{2}=(k+1)\phi^2(\Leftrightarrow h^{2}=(k+1)\phi^2).\label{2.6}
\eea
In \cite{dp}, Dileo and Pastore introduced the notion of $(k,\mu)'$-nullity distribution, on an almost $(2n+1)$-dimensional Kenmotsu manifold  $(M, \phi, \xi, \eta, g)$, which is defined for any $p \in M$ and $k,\mu \in \mathbb {R}$ as follows:
\bea
 N_{p}(k,\mu)'=\{Z\in T_{p}(M):R(X,Y)Z&=&k[g(Y,Z)X-g(X,Z)Y]\nonumber\\&&+\mu[g(Y,Z)h'X-g(X,Z)h'Y]\}.\label{2.7}
 \eea
The above notion is called generalized nullity distributions when one allows $k, \mu$ to be smooth functions.\\
Let $X\in\mathcal D$ be the eigen vector of $h'$ corresponding to the eigen value $\alpha$. Then from (\ref{2.6}) it is clear that $\alpha^{2} = -(k+1)$, a constant. Therefore $k\leq -1$ and $\alpha = \pm\sqrt{-k-1}$. We denote by $[\alpha]'$ and $[-\alpha]'$ the corresponding eigen spaces related to the non-zero eigen value $\alpha$ and $-\alpha$ of $h'$, respectively. In \cite{dp}, it is proved that in a $(2n+1)$ dimensional $(k,\mu)'$-almost Kenmotsu manifold $M$ with $h' \neq 0$, $k < -1,\;\mu = -2$ and Spec$(h') = \{0,\alpha,-\alpha\}$, with $0$ as simple eigen value and $\alpha = \sqrt{-k-1}$.  From (\ref{2.7}), we have
\bea
R(X,Y)\xi=k[\eta(Y)X-\eta(X)Y]+\mu[\eta(Y)h'X-\eta(X)h'Y],\label{2.8}
\eea
where
$k,\mu\in\mathbb R.$ Also we get from (\ref{2.8})
\bea
R(\xi,X)Y = k[g(X,Y)\xi - \eta(Y)X] + \mu[g(h'X,Y)\xi - \eta(Y)h'X].\label{2.9}
 \eea

Using (\ref{2.3}), we have
\bea
(\nabla_X \eta)Y = g(X,Y) - \eta(X)\eta(Y) + g(h'X,Y). \label{2.10}
\eea
For further details on $(k,\mu)'$-almost Kenmotsu manifolds, we refer the reader to go through the references (\cite{ucd}, \cite{dp}, \cite{pv}).

\section{\textbf{$\ast$-Miao-Tam critical equation}}

In this section, we study the notion of $\ast$-Miao-Tam critical equation in the framework of $(k,\mu)'$-almost Kenmotsu manifolds. To prove the main theorem, we need the following lemmas:
 
\begin{lemma}(\cite{xinxin}) \label{l3.1}
On a $(k,\mu)'$-almost Kenmotsu manifold with $k < -1$, the $\ast$-Ricci tensor is given by
\bea
S^\ast (X,Y) = -(k + 2)(g(X,Y) - \eta(X)\eta(Y)) \label{3.1}
\eea
for any vector fields $X$, $Y$.
\end{lemma}
Note that, $S^\ast$ is symmetric here and hence the notion of $\ast$-Miao-Tam critical equation is well defined in this setting. 

\begin{lemma} \label{l3.2}
If the metric of an almost contact metric manifold $(M,g)$ of dimension $n \geq 3$ satisfies the $\ast$-Miao-Tam critical equation, then the curvature tensor $R$ can be expressed as 
\bea
\nonumber R(X,Y)D\lambda &=& (Xf)Y - (Yf)X + (X\lambda)Q^\ast Y - (Y\lambda)Q^\ast X \\ && + \lambda [(\nabla_X Q^\ast)Y - (\nabla_Y Q^\ast)X] \label{3.2}
\eea
for any vector fields $X$ and $Y$ on $M$, where $f = - \frac{r^\ast \lambda + n}{n - 1}$.
\end{lemma}
\begin{proof}
Tracing (\ref{1.3}), we obtain $\Delta_g \lambda =  - \frac{r^\ast \lambda + n}{n - 1} = f$.\\
Therefore, equation (\ref{1.3}) can be written as
\bea
\nabla_X D\lambda = (1 + f)X + \lambda Q^\ast X. \label{3.3}
\eea
Differentiating the above equation covariantly along any vector field $Y$, we obtain
\bea
\nabla_Y \nabla_X D\lambda = (Yf)X + (1 + f)\nabla_Y X + (Y\lambda)Q^\ast X + \lambda \nabla_Y Q^\ast X. \label{3.4}
\eea
Interchanging $X$ and $Y$ in (\ref{3.4}) yields 
\bea
\nabla_X \nabla_Y D\lambda = (Xf)Y + (1 + f)\nabla_X Y + (X\lambda)Q^\ast Y + \lambda \nabla_X Q^\ast Y. \label{3.5}
\eea
Also, from (\ref{3.3}), we have
\bea
\nabla_{[X,Y]} D\lambda = ( 1 + f)(\nabla_X Y - \nabla_Y X) + \lambda Q^\ast(\nabla_X Y - \nabla_Y X). \label{3.6}
\eea
It is well known that,
\bea
\nonumber R(X,Y)D\lambda = \nabla_X \nabla_Y D\lambda - \nabla_Y \nabla_X D\lambda - \nabla_{[X,Y]}D\lambda. 
\eea
We now complete the proof by substituting (\ref{3.4})-(\ref{3.6}) in the foregoing equation.
\end{proof}

\begin{lemma} \label{l3.3}
The $\ast$-Ricci operator $Q^\ast$ of a $(2n + 1)$-dimensional $(k,\mu)'$-almost Kenmotsu manifold $M$ satisfy the following relation:
\bea
\nonumber (\nabla_X Q^\ast)Y - (\nabla_Y Q^\ast)X &=& (k + 2)[\eta(Y)(X - \eta(X)\xi - \phi hX) \\ && - \eta(X)(Y - \eta(Y)\xi - \phi hY)] \label{3.7}
\eea
for any vector fields $X$, $Y$ on $M$.
\end{lemma}
\begin{proof}
From (\ref{3.1}), we have 
\bea
Q^\ast X = - (k + 2)(X - \eta(X)\xi). \label{3.8}
\eea
Differentiating (\ref{3.8}) covariantly along any vector field $Y$ and using (\ref{2.3}), we obtain
\bea
\nabla_Y Q^\ast X = -(k + 2)[\nabla_Y X - (\nabla_Y \eta(X))\xi - \eta(X)(Y - \eta(Y)\xi - \phi hY)]. \label{3.9}
\eea
Now, with the help of (\ref{3.8}), (\ref{3.9}) and (\ref{2.10}), we calculate
\bea
\nonumber (\nabla_Y Q^\ast)X &=& \nabla_Y Q^\ast X - Q^\ast(\nabla_Y X) \\ \nonumber &=& (k + 2)[(g(X,Y) - \eta(X)\eta(Y) + g(h'X,Y))\xi \\ && + \eta(X)(Y - \eta(Y)\xi - \phi hY)]. \label{3.10}
\eea
In a similar manner, we obtain
\bea
\nonumber (\nabla_X Q^\ast)Y &=& (k + 2)[(g(X,Y) - \eta(X)\eta(Y) + g(h'X,Y))\xi \\ && + \eta(Y)(X - \eta(X)\xi - \phi hX)]. \label{3.11}
\eea
Now, subtracting (\ref{3.10}) from (\ref{3.11}), we obtain our desired result.
\end{proof}

\begin{theorem} \label{t3.4}
 If the metric of a $(2n + 1)$-dimensional
 $(k,\mu)'$-almost Kenmotsu manifold $M$ satisfies the $\ast$-Miao-Tam critical equation, then the manifold $M$ is $\ast$-Ricci flat and locally isometric to $\mathbb{H}^{n+1}(-4)$ $\times$ $\mathbb{R}^n$.
\end{theorem}

\begin{proof}
Tracing (\ref{3.1}), we have $r^\ast = - 2n(k + 2)$.\\
Using lemma \ref{l3.2} and lemma \ref{l3.3}, for a $(2n + 1)$-dimensional $(k,\mu)'$-almost Kenmotsu manifold $M$, we can write
\bea
\nonumber R(X,Y)D\lambda &=& (Xf)Y - (Yf)X + (X\lambda)Q^\ast Y - (Y\lambda)Q^\ast X \\ \nonumber && + (k + 2)[\eta(Y)(X - \eta(X)\xi - \phi hX) - \eta(X)(Y - \eta(Y)\xi - \phi hY)],
\eea
where $ f = - \frac{2n + 1 - 2n(k + 2)\lambda}{2n}$.\\
 Replacing $X$ by $\xi$ in the above equation yields
\bea
\nonumber R(\xi,Y)D\lambda = (\xi f)Y - (Yf)\xi + (\xi \lambda)Q^\ast Y - (k + 2)(Y - \eta(Y)\xi - \phi hY).
\eea
Taking inner product of the foregoing equation with $X$, we have
\bea
\nonumber g(R(\xi,Y)D\lambda,X) &=& (\xi f)g(X,Y) - (Yf)\eta(X) + (\xi \lambda)S^\ast(X,Y) \\ && - (k + 2)[g(X,Y) - \eta(X)\eta(Y) + g(h'X,Y)]. \label{3.12}
\eea
Again, since $g(R(\xi,Y)D\lambda,X) = - g(R(\xi,Y)X,D\lambda)$, then using (\ref{2.9}), we obtain
\bea
\nonumber g(R(\xi,Y)D\lambda,X) &=& - kg(X,Y)(\xi \lambda) + k\eta(X)(Y\lambda) \\ && + 2g(h'X,Y)(\xi \lambda) - 2\eta(X)((h'Y)\lambda). \label{3.13}
\eea
Equating (\ref{3.12}) and (\ref{3.13}) and then antisymmetrizing yields
\bea
\nonumber (Xf)\eta(Y) - (Yf)\eta(X) &=& k\eta(X)(Y\lambda) - k\eta(Y)(X\lambda) \\ && -2\eta(X)((h'Y)\lambda) + 2\eta(Y)((h'X)\lambda). \label{3.14}
\eea
Since $ f = - \frac{2n + 1 - 2n(k + 2)\lambda}{2n}$, then $(Xf) = (k + 2)(X\lambda)$ for any vector field $X$. Hence, equation (\ref{3.14}) reduces to
\bea
2(k + 1)[\eta(Y)(X\lambda) - \eta(X)(Y\lambda)] + 2\eta(X)((h'Y)\lambda) - 2\eta(Y)((h'X)\lambda) = 0. \label{3.15}
\eea
Setting $Y = \xi$ in (\ref{3.15}), we infer that
\bea
\nonumber (k + 1)[(X\lambda) - (\xi \lambda)\eta(X)] - ((h'X)\lambda) = 0,
\eea
which implies
\bea
(k + 1)[D\lambda - (\xi \lambda)\xi] - h'(D\lambda) = 0. \label{3.16}
\eea
Operating $h'$ on the above equation and using (\ref{2.6}), we obtain
\bea
h'(D\lambda) = - [ D\lambda - (\xi \lambda)\xi]. \label{3.17}
\eea
Substituting (\ref{3.17}) in (\ref{3.16}) yields
\bea
\nonumber (k + 2)[D\lambda - (\xi \lambda)\xi] = 0,
\eea
which implies either $ k = - 2$ or $ D\lambda = (\xi \lambda)\xi$.\\

\par Case 1: : If $k = -2$, then using $\alpha^2 = - k - 1$, we get $\alpha = \pm 1$. Without loss of generality, we consider $\alpha = 1$. Then from Prop. 4.2 of \cite{dp}, we have
\bea
\nonumber 	R(X_\alpha,Y_\alpha)Z_\alpha = - 4[g(Y_{\alpha},Z_{\alpha})X_{\alpha} - g(X_{\alpha},Z_{\alpha})Y_{\alpha}]
\eea
and
 \bea
 \nonumber 	R(X_{-\alpha},Y_{-\alpha})Z_{-\alpha} = 0
  \eea
  for any $X_\alpha,Y_\alpha,Z_\alpha \in [\alpha]'$ and $X_{-\alpha},Y_{-\alpha},Z_{-\alpha} \in [-\alpha]'$.  It is proved in Prop. 10 and Theorem 5 of \cite{dp2} that the distributions $[\xi] \oplus [\alpha]'$ and $[-\alpha]'$ are integrable and totally geodesic. Then it follows that the $(2n + 1)$-dimensional $(k,\mu)'$-almost Kenmotsu manifold is locally isometric to the direct product of an $(n + 1)$-dimensional integral manifold of $[\xi] \oplus [\alpha]'$ and an $n$-dimensional integral manifold of $[-\alpha]'$. Now, the restriction of the curvature tensor $R$ on $[\xi] \oplus [\alpha]'$ shows that the maximal integral manifold of $[\xi] \oplus [\alpha]'$ is a hyperbolic space of constant curvature $-4$. Also, the restriction of the curvature tensor $R$ on $[-\alpha]'$ shows that the maximal integral manifold of $[-\alpha]'$ is flat. Hence, we can say that the $(2n + 1)$-dimensional $(k,\mu)'$-almost Kenmotsu manifold $M$ is locally isometric to $\mathbb{H}^{n+1}(-4)$ $\times$ $\mathbb{R}^n$.\\
  Again from (\ref{3.1}), $k = -2$ implies $S^\ast = 0$ and therefore the manifold is $\ast$-Ricci flat.\\

\par Case 2: If $D\lambda = (\xi \lambda)\xi$, then differentiating this covariantly along any vector field $X$, we obtain
\bea
\nabla_X D\lambda = (X(\xi \lambda))\xi + (\xi \lambda)(X - \eta(X)\xi - \phi hX).\label{3.18}
\eea
Equating (\ref{3.3}) and (\ref{3.18}), we obtain
\bea
 Q^\ast X = \frac{1}{\lambda}((\xi \lambda) - 1 - f)X + \frac{1}{\lambda}(X(\xi \lambda) - (\xi \lambda)\eta(X))\xi - \frac{1}{\lambda}(\xi \lambda)\phi hX. \label{3.19}
\eea
Now, comparing the coefficients of $X$, $\xi$ and $\phi hX$ from (\ref{3.19}) and (\ref{3.8}), we get
\bea
(\xi \lambda) - 1 - f = - (k + 2)\lambda. \label{3.20}
\eea
\bea
X(\xi \lambda) - (\xi \lambda)\eta(X) = (k + 2)\lambda. \label{3.21}
\eea
\bea
(\xi \lambda) = 0. \label{3.22}
\eea
Using (\ref{3.22}) in (\ref{3.21}), we get $k = - 2$ as $\lambda$ is a non-zero smooth function. Therefore, from (\ref{3.20}), we have $ f = - 1$. Using $k = - 2$ in $f = - 1$, we will arrive at a contradiction.\\
This completes the proof.
\end{proof}

\section{\textbf{Example of a  $(k,\mu)'$-almost Kenmotsu manifold}}
In \cite{dp}, Dileo and Pastore  gives an example of a $(2n + 1)$-dimensional almost Kenmotsu manifold in which the characteristic vector field $\xi$ belongs to the $(k,\mu)'$-nullity distribution. For $5$-dimensional case, $k = -2$ (see \cite{dey2}). This implies that the manifold is locally isometric to $\mathbb{H}^3(-4) \times \mathbb{R}^2$ .\\
Again since $k = -2$, by (\ref{3.1}), the manifold is $\ast$-Ricci flat, that is, $S^\ast = 0$. Now tracing (\ref{1.3}), we obtain $\Delta_g \lambda = - \frac{5}{4}$. Therefore, $(g,\lambda)$ is a solution of the $\ast$-Miao-Tam critical equation, where $\lambda$ is given by the above Poisson equation. Hence, the theorem is verified.\\


\end{document}